\documentclass[10pt]{myelsart}
\usepackage{natbib}
\usepackage[leqno]{amsmath}
\usepackage{amsfonts}
\usepackage[all]{xy}
\usepackage{algorithm}
\usepackage{algorithmic}

\DeclareMathOperator{\depth}{depth}

\DeclareMathOperator{\der}{Der}
\DeclareMathOperator{\enm}{End}
\DeclareMathOperator{\ext}{Ext}
\DeclareMathOperator{\hc}{HC}
\DeclareMathOperator{\hh}{HH}
\DeclareMathOperator{\hmm}{Hom}
\DeclareMathOperator{\id}{id}
\DeclareMathOperator{\im}{im}

\DeclareMathOperator{\spec}{Spec}

\newcommand{\cc}{\Cset}
\newcommand{\cpd}{*}

\newcommand{\g}{\mathsf{g}}

\newcommand{\kdf}[1]{\Omega_{#1}}

\newcommand{\acl}{\mathsf{a}}
\newcommand{\ksker}{\mathsf{V}}
\newcommand{\lc}{\mathsf{lc}}
\newcommand{\vect}[1]{\mathbf{#1}}

\begin{document}

\begin{frontmatter}

\title{Computing obstructions for existence of connections on modules}
\author{Eivind Eriksen}
\address{Oslo University College}
\ead{eeriksen@hio.no}
\author{Trond St{\o}len Gustavsen}
\address{Buskerud University College}
\ead{trond.gustavsen@hibu.no}
\date{\today}

\begin{abstract}
We consider the notion of a connection on a module over a commutative
ring, and recall the obstruction calculus for such connections. This
obstruction calculus is defined using Hochschild cohomology. However,
in order to compute with Gr\"obner bases, we need the conversion to a
description using free resolutions. We describe our implementation in
{\sc Singular 3.0}, available as the library \texttt{conn.lib}.
Finally, we use the library to verify some known results and to obtain
a new theorem for maximal Cohen-Macaulay (MCM) modules on isolated
singularities. For a simple hypersurface singularity of dimension one
or two, it is known that all MCM modules admit connections. We prove
that for a simple threefold hypersurface singularity of type $A_n$,
$D_n$ or $E_n$, only the free MCM modules admit connections if $n \le
50$.
\end{abstract}

\end{frontmatter}

\section*{Introduction}

Let $k$ be an algebraically closed field of characteristic $0$. For
any commutative $k$-algebra $A$ and any $A$-module $M$, we
investigate whether there exist $A$-linear maps $\nabla: \g \to
\enm_k(M)$ satisfying the Leibniz rule
    \[ \nabla_D(am) = a \nabla_D(m) + D(a) \, m \; \text{ for all }
    D \in \g, \; a \in A, \; m \in M \]
on $\g = \der_k(A)$, or possibly on a smaller subset $\g \subseteq
\der_k(A)$ which is closed under the $A$-module and $k$-Lie algebra
structures of $\der_k(A)$. We refer to such maps as
$\g$-connections, or connections in case $\g = \der_k(A)$.

In many situations, this question is related to the topology of the
singularity $X = \spec(A)$. In fact, assume that $k$ is the field of
complex numbers, that $X$ is an isolated singularity and that $M$ is a
maximal Cohen-Macaulay (MCM) module. Then $M$ is locally free on the
complement $U \subset X$ of the singularity. If there is a connection
$\nabla: \der_k(A) \to \enm_k(M)$ that is also a Lie-algebra
homomorphism, we get an integrable connection on the vector bundle
$\widetilde{M}|_U$ on $U.$ Passing to the associated complex analytic
manifold $U^{an}$, we get a representation of the fundamental group
$\pi_{1}(U^{an})$ via the Riemann-Hilbert correspondence, see for
instance chapter 1 in \citet{del70} for the general case, or
\citet{gu-il06} for the case of normal surface singularities.

In this paper, we use algebraic methods to study the existence of
$\g$-connections on modules. To this aim, we recall the obstruction
calculus for $\g$-connections, which can be effectively implemented.
We use Hochschild cohomology to define this obstruction theory, but
for the implementation, the description of the obstructions in terms
of free resolutions, given in section \ref{s:freeres}, is essential.
We present our implementation as the library \texttt{conn.lib}
\citep{er-gu06-lib} for the computer algebra system {\sc Singular
3.0} \citep{gps05}.

In the case of simple hypersurface singularities (of type $A_n$, $D_n$
or $E_n$) in dimension $d$, there exists a connection on any MCM module
if $d \le 2$. Using our implementation, we get interesting results in
higher dimensions: For $d = 3$, we show that the only MCM modules that
admit connections are the free modules if $n \le 50$, and experimental
results indicate that the same result hold for $d = 4$.

These results led us to conjecture that for any simple hypersurface
singularity of dimension $d = 3$, the only MCM modules that admit
connections are the free modules. Using different techniques, we
proved a more general result in \citet{er-gu06}: An MCM module over
a simple hypersurface singularity of dimension $d \ge 3$ admits a
connection if and only if it is free.

In the case of simple elliptic surface singularities, it was shown in
\citet{ka88}, using analytic methods, that any MCM module admits a
connection. We verify some instances of this result
\emph{algebraically}, using our implementation.

\section{Basic definitions}

Let $k$ be an algebraically closed field of characteristic $0$, and let
$A$ be a commutative $k$-algebra. A \emph{Lie-Rinehart algebra} of
$A/k$ is a pair $(\g, \tau)$, where $\g$ is an $A$-module and a $k$-Lie
algebra, and $\tau: \g \to \der_k(A)$ is a morphism of $A$-modules and
$k$-Lie algebras, such that
    \[ [D, a D'] = a [D,D'] + \tau_D(a) D' \]
for all $D,D' \in \g$ and all $a \in A$, see \citet{ri63}. A
Lie-Rinehart algebra is the algebraic analogue of a \emph{Lie
algebroid}, and it is also known as a Lie pseudo-algebra or a
Lie-Cartan pair.

When $\g$ is a subset of $\der_k(A)$ and $\tau: \g \to \der_k(A)$ is
the inclusion map, the pair $(\g, \tau)$ is a Lie-Rinehart algebra if
and only if $\g$ is closed under the $A$-module and $k$-Lie algebra
structures of $\der_k(A)$. We are mainly interested in Lie-Rinehart
algebras of this type, and omit $\tau$ from the notation.

Let $\g$ be a Lie-Rinehart algebra. For any $A$-module $M$, we define a
\emph{$\g$-connection} on $M$ to be an $A$-linear map $\nabla: \g \to
\enm_k(M)$ such that
\begin{equation} \label{e:dp}
\nabla_D(am) = a \nabla_D(m) + D(a) m
\end{equation}
for all $D \in \g, \; a \in A, \; m \in M$. We say that $\nabla$
satisfies the \emph{derivation property} when condition (\ref{e:dp})
holds for all $D \in \g$. If $\nabla: \g \to \enm_k(M)$ is a $k$-linear
map that satisfies the derivation property, we call $\nabla$ a
$k$-linear $\g$-connection on $M$. A \emph{connection} on $M$ is a
$\g$-connection on $M$ with $\g = \der_k(A)$.

Let $\nabla$ be a $\g$-connection on $M$. We define the
\emph{curvature} of $\nabla$ to be the $A$-linear map $R_{\nabla}: \g
\wedge \g \to \enm_A(M)$ given by $R_{\nabla}(D \wedge D') = [
\nabla_D, \nabla_{D'} ] - \nabla_{[D,D']}$ for all $D,D' \in \g$. We
say that $\nabla$ is an \emph{integrable $\g$-connection} if
$R_{\nabla} = 0$.

When $A$ is a regular $k$-algebra, it is usual to define a
\emph{connection} on $M$ to be a $k$-linear map $\nabla: M \to M
\otimes_A \kdf A$ such that $\nabla(am) = a \nabla(m) + m \otimes d(a)$
for all $a \in A, \; m \in M$, see \citet{kat70}. Moreover, the
\emph{curvature} of $\nabla$ is usually defined as the $A$-linear map
$R_{\nabla}: M \to M \otimes_A \kdf{A}^2$ given by $R_{\nabla} =
\nabla^1 \circ \nabla$, where $\nabla^1$ is the natural extension of
$\nabla$ to $M \otimes \kdf A$, and $\nabla$ is an \emph{integrable
connection} if $R_\nabla = 0$.

Let $A$ be any commutative $k$-algebra. For expository purposes, we
define an \emph{$\Omega$-connection} on an $A$-module $M$ to be a
connection on $M$ in the sense of the preceding paragraph. By the
universal property of $\kdf A$, it follows that any (integrable)
$\Omega$-connection on $M$ induces an (integrable) connection on $M$.

When $A$ is a regular $k$-algebra essentially of finite type, there
is a bijective correspondence between connections on $M$ and
$\Omega$-connections on $M$ for any finitely generated $A$-module
$M$. However, there are many modules that admit connections but not
$\Omega$-connections when $A$ is a singular $k$-algebra, see
subsection \ref{ss:runtime} for some examples.

\section{Hochschild cohomology}

Let $k$ be an algebraically closed field of characteristic $0$, and let
$A$ be any commutative $k$-algebra. For any $A$-$A$ bimodule $Q$, the
\emph{Hochschild cohomology} $\hh^{\cpd}(A, Q)$ of $A$ with values in
$Q$ is the cohomology of the Hochschild complex $( \hc^{\cpd}, d^{\cpd}
)$, where
    \[ \hc^n(A,Q) = \hmm_k(\otimes^n_k A, Q) \]
We recall that $d^0: Q \to \hmm_k(A,Q)$ is given by $d^0(q)(a) = aq-qa$
for all $a \in A, \; q \in Q$ and that $d^1: \hmm_k(A,Q) \to \hmm_k(A
\otimes_k A,Q)$ is given by $d^1(\phi)(a \otimes b) = a \phi(b) -
\phi(ab) + \phi(a) b$ for all $a,b \in A, \; q \in Q$. For the
definitions of the higher differentials, we refer to \citet{wei94}. It
is clear that $\ker(d^1) = \der_k(A,Q)$, and we refer to $\im(d^0)$ as
the \emph{trivial derivations}.

Given $A$-modules $M, M'$, we shall consider Hochschild cohomology of
$A$ with values in the bimodule $Q = \hmm_k(M,M')$. In this case, the
natural exact sequence $0 \to \hh^0(A,Q) \to \hc^0(A,Q) \to \ker(d^1)
\to \hh^1(A,Q) \to 0$ has the form
\begin{multline} \label{s:hh}
0 \to \hmm_A(M,M') \to \hmm_k(M,M') \xrightarrow{d^0} \\
\to \der_k(A,\hmm_k(M,M')) \to \hh^1(A, \hmm_k(M,M')) \to 0
\end{multline}

It follows from \citet{wei94}, theorem 8.7.10 and lemma 9.1.9, that
there is a natural isomorphism of $k$-linear vector spaces
\begin{equation} \label{e:extid}
\ext^n_A(M,M') \to \hh^n(A,\hmm_k(M,M'))
\end{equation}
for any $n \ge 0$. For our purposes, it is useful to describe an
identification of this type explicitly in the case $n = 1$.

Let $(L_{\cpd}, d_{\cpd})$ be a free resolution of $M$, and let $\tau:
M \to L_0$ be a $k$-linear section of the augmentation morphism $\rho:
L_0 \to M$. Let moreover $\phi: L_1 \to M'$ be a $1$-cocycle in
$\hmm_A(L_{\cpd}, M')$. For any $a \in A, \; m \in M$, we can find an
element $x = x(a,m) \in L_1$ such that $d_0(x) = a \tau(m) - \tau(am)$.
It is easy to check that $\psi(a)(m) = \phi(x)$ defines a derivation
$\psi \in \der_k(A, \hmm_k(M,M'))$ which is independent of the choice
of $x$.

\begin{lem} \label{l:id}
The assignment $\phi \mapsto \psi$ defined above induces a natural
injective map $\sigma: \ext^1_A(M,M') \to \hh^1(A, \hmm_k(M,M'))$ of
$k$-linear vector spaces.
\end{lem}
\begin{pf}
If $\phi = \phi' d_0$ is a coboundary in $\hmm_A(L_{\cpd},M')$, then
$\psi$ is a trivial derivation, given by $\psi(a)(m) = a \psi'(m) -
\psi'(am)$ with $\psi' = \phi' \tau$. Hence the assignment induces a
well-defined map $\sigma$ of $k$-linear vector spaces. To see that it
is injective, assume that $\phi$ is a $1$-cocycle in $\hmm_A(L_{\cpd},
M')$ which maps to a trivial derivation $\psi$, given by $\psi(a)(m) =
a \psi'(m) - \psi'(am)$ for some $\psi' \in \hmm_k(M,M')$. For any $x
\in L_0$, there exists an element $x' \in L_1$ such that $d_0(x') = x -
\tau \rho(x)$, and $x'$ is unique modulo $\im(d_1)$. Let us define
$\phi': L_0 \to M'$ to be the map given by $\phi'(x) = \psi' \rho(x) +
\phi(x')$, this is clearly a well-defined map since $\phi \, d_1 = 0$.
One may show that $\phi'$ is $A$-linear, and it satisfies $\phi' d_0 =
\phi$ by construction. Hence $\phi$ is a coboundary.
\end{pf}

We do not claim that $\sigma$ coincides with the identification
(\ref{e:extid}) for $n = 1$. We shall use $\sigma$ to obtain a concrete
identification of certain classes defined using free resolutions, see
section \ref{s:freeres} for details, with the obstruction classes that
we define in the next section. The important fact is therefore that
$\sigma$ is injective.

\section{Obstruction theory}
\label{s:obstr}

Let $k$ be an algebraically closed field of characteristic $0$, let $A$
be a commutative $k$-algebra, and let $M$ be an $A$-module.

\begin{prop} \label{p:acl}
There is a canonical class $\acl(M) \in \ext^1_A(M,M \otimes_A \kdf
A)$, called the Atiyah class of $M$, such that $\acl(M) = 0$ if and
only if there is an $\Omega$-connection on $M$. In this case, there is
a transitive and effective action of $\hmm_A(M,M \otimes_A \kdf A)$ on
the set of $\Omega$-connections on $M$.
\end{prop}
\begin{pf}
We consider the derivation $\psi: A \to \hmm_k(M,M \otimes_A \kdf A)$
given by $\psi(a)(m) = m \otimes d(a)$, where $d$ is the universal
derivation of $A$, and define $\acl(M)$ to be the class in
$\ext^1_A(M,M \otimes_A \kdf A)$ corresponding to the class $[\psi]$ of
$\psi$ in $\hh^1(A,\hmm_k(M,M \otimes_A \kdf A))$ via the
identification (\ref{e:extid}). Using the sequence (\ref{s:hh}), the
proposition follows easily.
\end{pf}

\begin{prop} \label{p:ksker}
There is a canonical map $g: \der_k(A) \to \ext^1_A(M,M)$, called the
Kodaira-Spencer map of $M$, with the following properties:
\begin{enumerate}
\item The Kodaira-Spencer kernel $\ksker(M) = \ker(g)$ is
a Lie algebroid of $A/k$,
\item For any $D \in \der_k(A)$, there exists an operator $\nabla_D
\in \enm_k(M)$ with derivation property with respect to $D$ if and only
if $D \in \ksker(M)$.
\end{enumerate}
In particular, $\ksker(M)$ is maximal among the subsets $\g \subseteq
\der_k(A)$ such that there exists a $k$-linear $\g$-connection on $M$.
\end{prop}
\begin{pf}
For any $D \in \der_k(A)$, we consider $\psi_D \in \der_k(A,\enm_k(M))$
given by $\psi_D(a)(m) = D(a) m$, and denote by $g(D)$ the class in
$\ext^1_A(M,M)$ corresponding to the class $[\psi_D]$ of $\psi_D$ in
$\hh^1(A, \enm_k(M))$ via the identification (\ref{e:extid}). This
defines the Kodaira-Spencer map $g$ of $M$, which is $A$-linear by
definition. Clearly, its kernel $\ksker(M)$ is closed under the Lie
product. Using the exact sequence (\ref{s:hh}), it easily follows that
there exists an operator $\nabla_D$ with derivation property with
respect to $D$ if and only if $D \in \ksker(M)$.
\end{pf}

We remark that the Kodaira-Spencer map $g: \der_k(A) \to \ext^1_A(M,M)$
is the contraction against the Atiyah class $\acl(M) \in \ext^1_A(M,M
\otimes_A \kdf A)$. See also \citet{kal05}, section 2.2 for another
proof of proposition \ref{p:ksker}.

\begin{prop} \label{p:lc}
There is a canonical class $\lc(M) \in \ext^1_A(\ksker(M), \enm_A(M))$
such that $\lc(M) = 0$ if and only if there exists a
$\ksker(M)$-connection on $M$. In this case, there is a transitive and
effective action of $\hmm_A(\ksker(M), \enm_A(M))$ on the set of
$\ksker(M)$-connections on $M$.
\end{prop}
\begin{pf}
Let $\ksker = \ksker(M)$, choose a $k$-linear $\ksker$-connection
$\nabla: \ksker \to \enm_k(M)$ on $M$, and let $\phi \in \der_k(A,
\hmm_k(\ksker, \enm_A(M)))$ be the derivation given by $\phi(a)(D) = a
\nabla_{D} - \nabla_{aD}$. We denote by $\lc(M)$ the class in
$\ext^1_A(\ksker, \enm_A(M))$ corresponding to the class $[\phi]$ of
$\phi$ in $\hh^1(A, \hmm_k(\ksker, \enm_A(M)))$ via the identification
(\ref{e:extid}). One may check that this class is independent of
$\nabla$. Using the exact sequence (\ref{s:hh}), the proposition
follows easily.
\end{pf}

There is a natural short exact sequence $0 \to \enm_A(M) \to c(M) \to
\ksker(M) \to 0$ of left $A$-modules, where $c(M) = \{ P \in \enm_k(M):
[P,a] \in A \text{ for all } a \in A \}$ is the module of first order
differential operators on $M$ with scalar symbol, and $c(M) \to
\der_k(A)$ is the natural map, given by $P \mapsto [P,-]$, with image
$\ksker(M)$, see also \citet{kal05}, proposition 2.2.10. We remark that
this extension of left $A$-modules splits if and only if $\lc(M) = 0$.

\begin{lem} \label{l:dirsum}
For any $A$-modules $M,M'$, we have $\ksker(M \oplus M') = \ksker(M)
\cap \ksker(M')$, and $\lc(M \oplus M') = 0$ if and only if $\lc(M) =
\lc(M') = 0$.
\end{lem}

Lemma \ref{l:dirsum} is a direct consequence of proposition
\ref{p:ksker} and \ref{p:lc}. We remark that the first part also
follows from \citet{bu-li04}, lemma 3.4.

\section{Computing with free resolutions}
\label{s:freeres}

We define new classes $\acl(M)' \in \ext^1_A(M,M \otimes_A \kdf A)$,
$g'(D) \in \ext^1_A(M,M)$ for all $D \in \der_k(A)$, and $\lc(M)' \in
\ext^1_A(\ksker(M),\enm_A(M))$ in this section, and show that via
$\sigma$, these classes correspond to the obstructions $\acl(M)$,
$g(D)$ and $\lc(M)$ defined in section \ref{s:obstr}. Since the new
classes are defined using free resolutions, we can use Gr\"obner bases
to compute them in {\sc Singular}.

Let $k$ be an algebraically closed field of characteristic $0$, let $A$
be a commutative $k$-algebra essentially of finite type, and let $M$ be
a finitely generated $A$-module. Let
    \[ 0 \gets M \xleftarrow{\rho} L_0 \xleftarrow{d_0} L_1
    \gets \dots \]
be a free resolution of $M$ such that $L_i$ has finite rank for all $i
\ge 0$. We choose bases $\{ e_1, \dots, e_m \}$ of $L_0$ and $\{ f_1,
\dots, f_n \}$ of $L_1$, and write $( a_{ij} )$ for the matrix of $d_0:
L_1 \to L_0$ with respect to the chosen bases. One may show that the
matrix $( d(a_{ij}) )$, considered as an $A$-linear map $L_1 \to L_0
\otimes_A \kdf A$, defines a $1$-cocycle in $\hmm_A(L_{\cpd},M
\otimes_A \kdf A)$ and therefore a cohomology class $\acl(M)' \in
\ext^1_A(M,M \otimes_A \kdf A)$. The class $\acl(M)'$ is called the
\emph{Atiyah class} of $M$, see \citet{an-lj89}.

For any $D \in \der_k(A)$, the matrix $( D(a_{ij}) )$, considered as an
$A$-linear map $L_1 \to L_0$, defines a $1$-cocycle in
$\hmm_A(L_{\cpd},M)$ and therefore a cohomology class $g'(D) \in
\ext^1_A(M,M)$. We remark that the map $D \mapsto g'(D)$ is the
contraction against $\acl(M)'$, and it is therefore well-known that
$g'$ is the Kodaira-Spencer map of $M$, see \citet{ill71}.

Since $A$ is essentially of finite type over $k$, it follows that
$\ksker(M)$ is a left $A$-module of finite type. Let
    \[ 0 \gets \ksker(M) \xleftarrow{\rho'} L'_0 \xleftarrow{d'_0}
    L'_1 \gets \dots \]
be a free resolution of $\ksker(M)$ such that $L'_i$ has finite rank
for all $i \ge 0$. We choose bases $\{ e'_1, \dots, e'_p \}$ of $L'_0$
and $\{ f'_1, \dots, f'_q \}$ of $L'_1$, and write $( c_{ij} )$ for the
matrix of $d'_0: L'_1 \to L'_0$ with respect to the chosen bases. For
$1 \le i \le p$, let us write $D_i = \rho'(e'_i)$. Since $D_i \in
\ksker(M)$, we can find a $k$-linear operator $\nabla_i: M \to M$ which
has the derivation property with respect to $D_i$ by proposition
\ref{p:ksker}. We allow a slight abuse of notation, and write $D_i: L_0
\to L_0$ for the action given by $D_i(a e_j) = D(a) e_j$ for all $a \in
A$ and $1 \le j \le n$. Since the map
    \[ \nabla_i \, \rho - \rho \; D_i: L_0 \to M \]
is $A$-linear, we can lift it to an $A$-linear map $P_i: L_0 \to L_0$
via $\rho$, i.e. such that $\nabla_i \, \rho - \rho D_i = \rho P_i$,
and this implies that $\nabla_i \, \rho = \rho (D_i + P_i)$.

Let us consider $R_j = \sum_i c_{ij} P_i$ as an $A$-linear map $R_j:
L_0 \to L_0$ for $1 \le j \le q$. Since $\sum_i c_{ij} D_i = 0$ for all
$j$, we see that $R_j$ induces an $A$-linear operator on $M$ for $1 \le
j \le q$. One may show that the $A$-linear map $L'_1 \to \enm_A(L_0)$
given by $f'_j \mapsto R_j$ induces a $1$-cocycle in $\hmm_A(L'_{\cpd},
\enm_A(M))$, and therefore defines a cohomology class $\lc(M)'$ in
$\ext^1_A(\ksker(M), \enm_A(M))$.

\begin{lem} \label{l:idobstr}
Let $\acl(M)$, $g(D)$ and $\lc(M)$ be the obstructions defined in
section \ref{s:obstr}. Then we have:
\begin{enumerate}
\item $\acl(M) = 0$ if and only if $\acl(M)' = 0$,
\item $\ker(g) = \ker(g')$,
\item $\lc(M) = 0$ if and only if $\lc(M)'=0$.
\end{enumerate}
\end{lem}
\begin{pf}
To prove the first two parts, it is enough to show that $\acl(M)$ and
$\acl(M)'$ maps to the same element in $\hh^1(A, \hmm_k(M,M \otimes_A
\kdf A))$, since $\sigma$ is functorial, $g$ is the contraction against
$\acl(M)$, and $g'$ is the contraction against $\acl(M)'$. Since
$\acl(M)$ maps to $[ \psi ] \in \hh^1(A, \hmm_k(M,M \otimes_A \kdf
A))$, the class of the derivation $\psi$ defined in the proof of
proposition \ref{p:acl}, we must show that $\sigma(\acl(M)') = [\psi]$.
This follows from the fact that $( d(a_{ij}) ) = d d_0 - d_0 d$, where
$d: L_i \to L_i \otimes_A \kdf A$ is the natural action of the
universal derivation $d$ on $L_i$ for $i = 0,1$ with respect to the
chosen bases.

To prove the third part, is is enough to show that $\lc(M)$ and
$\lc(M)'$ maps to the same element in
$\hh^1(A,\hmm_k(\ksker(M),\enm_A(M)))$. Since $\lc(M)$ maps to $[ \phi
] \in \hh^1(A, \hmm_k(\ksker(M),\enm_A(M)))$, the class of the
derivation $\phi$ defined in the proof of proposition \ref{p:lc}, we
must show that $\sigma(\lc(M)') = [ \phi ]$. Let us consider the
$A$-linear map $\nabla: L'_0 \to \enm_k(M)$ given by $e'_i \mapsto
\nabla_i$ for $1 \le i \le p$, and notice that $\nabla d'_0(f'_j) =
R_j$ for $1 \le j \le q$ and that $\nabla \tau'$ is a $k$-linear
$\ksker(M)$-connection on $M$ for any $k$-linear section $\tau'$ of
$\rho'$. For any $a \in A, \; D \in \ksker(M)$, choose $x \in L'_1$
such that $d'_0(x) = a \tau'(D) - \tau'(aD)$. Then the cocycle $f'_j
\mapsto R_j$ maps $x$ to $\nabla(a \tau'(D) - \tau'(aD)) = a \nabla
\tau'(D) - \nabla \tau'(aD)$. By the definition of $\sigma$, it follows
that $\sigma(\lc(M)') = [ \phi ]$.
\end{pf}

We remark that the first and second part of lemma \ref{l:idobstr}
follows from more general results, see for instance \citet{bu-fl03}.

\section{Implementation} \label{s:algo}

Let $k$ be an algebraically closed field of characteristic $0$, let
$A$ be a $k$-algebra essentially of finite type, and let $M$ be a
finitely generated $A$-module. In this section, we explain the
procedures \texttt{Der()}, \texttt{AClass(M)}, \texttt{KSKernel(M)}
and \texttt{LClass(M)} in the {\sc Singular} library
\texttt{conn.lib} \citep{er-gu06-lib}, which calculate the
derivation module $\der_k(A)$, the Atiyah class $\acl(M)'$, the
Kodaira-Spencer kernel $\ksker(M)$ and the class $\lc(M)'$. We use
the notation from section \ref{s:freeres}.

In the implementation, we identify elements of free $A$-modules (with
given bases) with column vectors, and identify $A$-linear maps of free
$A$-modules (with given bases) with left multiplication by matrices. We
use presentation matrices to represent modules, and assume that $A$ is
the base ring in {\sc Singular}.

\texttt{Der()}: Let us recall how to compute $\der_k(A)$. If $A$ is
of finite type over $k$, we may assume that $A = S/I$, where $S =
k[x_1, \dots, x_n]$ and $I = (F_1, \dots, F_m) \subseteq S$ is an
ideal. The Jacobian matrix $J = ( \partial F_i / \partial x_j )$
defines an $A$-linear map $J: A^n \to A^m$, and it is well-known
that $\vect a \mapsto a_1 \partial / \partial x_1 + \dots + a_n
\partial / \partial x_n$ defines an isomorphism $\ker(J) \to
\der_k(A)$ of left $A$-modules. Finally, if $B$ is a localization of
$A$, then $\der_k(B) \cong B \otimes_A \der_k(A)$, so $\der_k(B)
\cong \ker(B \otimes_A J)$. In the implementation, we write
\texttt{j} for the matrix $J$ (respectively $B \otimes_A J$), and
\texttt{d} for the matrix with column space $\ker(J)$ (respectively
$\ker(B \otimes_A J)$). \texttt{Der()} returns the matrix
$\mathtt{d}$.

\texttt{AClass(M)}: In section \ref{s:freeres}, we have seen that
$\acl(M)' \in \ext^1_A(M, M \otimes_A \kdf A)$ is represented by $(
d(a_{ij}) ): L_1 \to L_0 \otimes_A \kdf A$, where $d: A \to \kdf A$ is
the universal derivation, and $( a_{ij} )$ is the matrix of $d_0$ with
respect to chosen bases of $L_0$ and $L_1$. Choose a free resolution
    \[ 0 \gets \kdf A \xleftarrow{\rho''} L''_0 \xleftarrow{d''_0} L''_1
    \gets \dots \]
of $\kdf A$, and lift $( d(a_{ij}) )$ to an $A$-linear map $\Delta: L_1
\to L_0 \otimes_A L''_0$. It follows that $\acl(M)'=0$ if and only if
$\Delta \in H_0$, where $H_0 \subseteq \hmm_A(L_1,L_0 \otimes_A L''_0)$
is the submodule given by $H_0 = \im(\hmm_A(d_0,L_0 \otimes_A L''_0)) +
\im(\hmm_A(L_1, \delta_0))$ and $\delta_0 = d_0 \otimes \id - \id
\otimes d''_0: (L_1 \otimes_A L''_0) \oplus (L_0 \otimes_A L''_1) \to
L_0 \otimes_A L''_0$ is the differential in the complex $L_\cpd
\otimes_A L''_\cpd$. In the implementation, we write \texttt{kPres} for
the matrix of $\delta_0$, \texttt{rel} for the matrix with column space
$H_0$, and \texttt{cm} for the representative of $\Delta$ with respect
to suitable bases. \texttt{AClass(M)} returns $( \mathtt{cm} \not \in
\mathtt{rel} )$.

\texttt{KSKernel(M)}: For any $D \in \der_k(A)$, we have seen that
$g'(D) \in \ext^1_A(M,M)$ is represented by $( D(a_{ij}) ): L_1 \to
L_0$, where $( a_{ij} )$ is the matrix of $d_0$ with respect to chosen
bases of $L_0$ and $L_1$. We write $\Delta: \der_k(A) \to
\hmm_A(L_1,L_0)$ for the $A$-linear map given by $D \mapsto
(D(a_{ij}))$. It follows that $g'(D) = 0$ if and only if $\Delta(D) \in
H_0$, where $H_0 \subseteq \hmm_A(L_1,L_0)$ is the submodule given by
$H_0 = \im(\hmm_A(d_0,L_0)) + \im(\hmm_A(L_1, d_0))$. Choose a free
resolution
    \[ 0 \gets \der_k(A) \xleftarrow{\rho'''} L'''_0 \xleftarrow{d'''_0}
    L'''_1 \gets \dots \]
of $\der_k(A)$, and lift $\Delta$ to an $A$-linear map $\gamma: L'''_0
\to \hmm_A(L_1,L_0)$. It follows that $\ksker(M) =
\rho'''(\gamma^{-1}(H_0)) \subseteq \der_k(A)$. In the implementation,
we write \texttt{rel} for the matrix with column space $H_0$,
\texttt{g} for the matrix of $\gamma$ and $\texttt{ker}$ for the matrix
with column space $\gamma^{-1}(H_0)$ with respect to suitable bases.
\texttt{KSKernel(M)} returns $( \mathtt{Der()} \cdot \mathtt{ker} \not
= \mathtt{Der()} )$.

\texttt{LClass(M)}: There is a $k$-linear connection $\nabla: \ksker
\to \enm_k(M)$ on $\ksker = \ksker(M)$ by proposition \ref{p:ksker},
and $\nabla_{D_i}$ can be lifted to a $k$-linear operator $D_i + P_i:
L_0 \to L_0$ for some $P_i \in \enm_A(L_0)$ for $1 \le i \le p$, where
we write $D_i = \rho'(e'_i)$. Let $\Lambda: L'_0 \to \hmm_A(L_0,L_0)$
be the $A$-linear map given by $e'_i \mapsto P_i$. Then $\Lambda d'_0$
induces an $A$-linear map $L'_1 \to \enm_A(M)$, since any $\phi \in
\im(\Lambda d'_0)$ satisfies $\phi d_0 = d_0 \psi$ for some $\psi \in
\hmm_A(L_1,L_1)$, and we have seen that $\Lambda d'_0$ represents
$\lc(M)' \in \ext^1_A(\ksker,\enm_A(M))$. Choose a free resolution
    \[ 0 \gets \enm_A(M) \xleftarrow{\rho''''} L''''_0
    \xleftarrow{d''''_0} L''''_1 \gets \dots \]
of $\enm_A(M)$, and lift $\Lambda d'_0$ to an $A$-linear map $\lambda:
L'_1 \to L''''_0$. It follows that $\lc(M)'=0$ if and only if $\lambda
\in H_0$, where $H_0  \subseteq \hmm_A(L'_1,L'''_0)$ is the submodule
$H_0 = \im(\hmm_A(d'_0,L'''_0)) + \im(\hmm_A(L'_1, d'''_0))$. In the
implementation, we write \texttt{vRel} for the matrix with column space
$H_0$, \texttt{L} for the matrix of $\Lambda$, and \texttt{lc} for the
representative of $\lambda$ with respect to suitable bases.
\texttt{LClass(M)} returns $( \mathtt{lc} \not \in \mathtt{vRel} )$.

\section{Results on maximal Cohen-Macaulay modules}

Let $k$ be an algebraically closed field of characteristic $0$, and let
$(A,m)$ be a local complete commutative Noetherian $k$-algebra. We say
that a finitely generated $A$-module $M$ is a maximal Cohen-Macaulay
(MCM) module if $\depth(M) = \dim(A)$. In this section, we give some
results on the existence of connections on MCM modules over isolated
singularities.

The {\sc Singular} library \texttt{conn.lib} \citep{er-gu06-lib} can
be used to calculate the obstructions $\ksker(M)$ and $\lc(M)$ for
the existence of connections on an $A$-module $M$ when $A$ is a
$k$-algebra of finite type or a local $k$-algebra essentially of
finite type, as explained in section \ref{s:algo}, using the
monomial ordering \texttt{dp} (or \texttt{lp}) respectively
\texttt{ds}. However, there is no direct way of computing with
complete algebras in {\sc Singular}. Nevertheless, thanks to the
following result, all examples in this section can be verified using
the library \texttt{conn.lib}:

\begin{lem}
Let $A$ be a commutative $k$-algebra of finite type with a
distinguished maximal ideal $m$, and let $M$ be a finitely generated
$A$-module. Consider $M_m$ as a module over the local ring $A_m$ and
$\widehat M$ as a module over the $m$-adic completion $\widehat A$.
If $M$ is locally free on $\spec A \setminus \{ m \}$, then
\begin{enumerate}
\item $\ksker(M) = \der_k(A) \Leftrightarrow \ksker(M_m) =
\der_k(A_m) \Leftrightarrow \ksker(\widehat M) = \der_k(\widehat
A)$,
\item $\lc(M) = 0 \Leftrightarrow \lc(M_m) = 0 \Leftrightarrow
\lc(\widehat M) = 0$.
\end{enumerate}
\end{lem}

\subsection{Simple hypersurface singularities}

Let $A$ be the complete local ring of a simple hypersurface
singularity, see \citet{Arn81} and \citet{Wal84}, and let $d \ge 1$
be the dimension of $A$. Then $A \cong k[[x, y, z_1, \dots,
z_{d-1}]] / (f)$, where $f$ is of the form
\begin{align*}
A_n: \; & f = x^2 + y^{n+1} + z_1^2 + \dots + z_{d-1}^2 & n \ge 1\\
D_n: \; & f = x^2y + y^{n-1} + z_1^2 + \dots + z_{d-1}^2 & n \ge 4\\
E_6: \; & f = x^3 + y^4 + z_1^2 + \dots + z_{d-1}^2 & \\
E_7: \; & f = x^3 + xy^3 + z_1^2 + \dots + z_{d-1}^2 & \\
E_8: \; & f = x^3 + y^5 + z_1^2 + \dots + z_{d-1}^2 &
\end{align*}
It is known that the simple hypersurface singularities are exactly the
hypersurface singularities which are of finite CM representation type,
i.e. which have a finite number of isomorphism classes of
indecomposable MCM $A$-modules, see \citet{Kno87} and
\citet{BucGreSch87}. The MCM modules are given by matrix
factorizations, see \citet{Eis80}, and are completely classified. A
complete list can be found in \citet{gr-kn85} in the curve case and
obtained from the McKay correspondence in the surface case. In higher
dimensions, this list can be obtained using Kn\"orrer's periodicity
theorem, see \citet{Kno87} and \citet{Sch87}. See also \citet{yos90}
for an overview.

We mention the following result for completeness. In the curve case,
the result appears in \citet{er06} and all computations were originally
done by hand constructing a connection on each indecomposable MCM
module. Now we can check the result for $E_6$, $E_7$, $E_8$ and for
$A_n$, $D_n$ for $n \le 50$ using {\sc Singular}. In the surfaces case,
the result is known to specialists, but has not been published as far
as we know. The simple hypersurface singularities are the quotient
surface singularities that are hypersurfaces. Then the result follows
from an observation of J. Christophersen, see \citet{er-gu06} or
\citet{maa05} for details.

\begin{thm}
Let $A$ be the complete local ring of a simple hypersurface singularity
(of type $A_n$, $D_n$ or $E_n$) of dimension $d \le 2$. Then there is a
connection on any MCM $A$-module.
\end{thm}

In contrast to the situation in dimension one and two, experimental
data suggest that there are very few MCM modules over simple
singularities of dimension $d \ge 3$ that admit connections. We have
discovered the following result using {\sc Singular}:

\begin{thm}
Let $A$ be the complete local ring of a simple hypersurface
singu\-larity (of type $A_n$, $D_n$ or $E_n$) of dimension $d = 3$, and
let $M$ be an MCM $A$-module. If $n \le 50$, then there is a connection
on $M$ if and only if $M$ is free.
\end{thm}

Experimental results indicate that the same result hold in dimension $d
= 4$. Due to the computational cost, we have not been able to calculate
$\lc(M)$ for any non-free MCM module $M$ over simple singularities in
dimension $d \ge 5$ using our implementation.

These results led us to conjecture that for any simple hypersurface
singularity of dimension $d = 3$, the only MCM modules that admit
connections are the free modules. Using different techniques, we
proved a more general result in \citet{er-gu06}: An MCM module over
a simple hypersurface singularity of dimension $d \ge 3$ admits a
connection if and only if it is free.

\subsection{Simple elliptic surface singularities}

Let $k = \cc$ and let $A$ be the complete local ring of a simple
elliptic surface singularity over $k$. Then $A$ is of tame CM
representation type, and all indecomposable MCM $A$-modules have been
classified, see \citet{ka89}. Moreover, it follows from the results in
\citet{ka88} that any MCM $A$-module admits a connection. Kahn obtained
these results by analytic methods.

We remark that we can prove many instances of Kahn's result
\emph{algebraically} using the library \texttt{conn.lib}. For instance,
when $A = k[[x,y,z]]/(x^3+y^3+z^3)$, we can use presentation matrices
of MCM $A$-modules from \citet{LazPfiPop02}.

\subsection{Other experimental results and timings}
\label{ss:runtime}

In the table below, we list the obstructions for existence of
connections in some additional cases. We used {\sc Singular 3.0} on
Cygwin in a Windows XP Professional (Version 2002, Service Pack 2)
environment, running on a PC (AMD Athlon 64 3800+ CPU 2.41 GHz, 1 GB
RAM). In one case, we were not able to calculate the obstructions
due to the computational cost. The timings are the total runtimes
for computations of all three obstructions.

The singularities in the table are: (1) The threefold non-Gorenstein
scroll of type $(2,1)$, see \citet{yos90}, proposition 16.12, (2) The
threefold non-Gorenstein quotient singularity of finite CM
representation type, see \citet{yos90}, proposition 16.10, (3) The
cubic cone with equation $x^3 + y^3 + z^3 = 0$, see
\citet{LazPfiPop02}, (4) The non-Gorenstein monomial curve with
numerical semigroup $(3,4,5)$, and (5) The threefold $E_6$ singularity.
The modules in the table are given in the references \citet{yos90} and
\citet{LazPfiPop02}, pp. 234--235, in the first $8$ cases, and as the
ideals $p = (t^3,t^4) \subseteq k[[t^3,t^4,t^5]]$ and $m = (x,y,z,w)
\subseteq k[[x,y,z,w]]/(x^3+y^4+z^2+w^2)$ in the last two cases.

We remark that the module $m$ is torsion free of depth two, while all
other modules in the table are MCM, and make the following
observations:

\begin{enumerate}
\item There is an isolated curve singularity with complete local
ring $A$ and an MCM $A$-module $M$ such that $M$ \emph{does not} admit
a connection.
\item We have not found an isolated surface singularity with
complete local ring $A$ and an MCM $A$-module $M$ such that $M$
\emph{does not} admit a connection.
\item There is an isolated threefold singularity with complete local
ring $A$ and a non-free MCM $A$-module $M$ such that $M$ admits a
connection.
\item It does not seem difficult to find an isolated singularity with
complete local ring $A$ and a (non-free) torsion free $A$-module $M$
such that $M$ admits a connection.
\end{enumerate}

\begin{center}
\begin{tabular}{|l|l|l|l|l|l|}
\hline
Singularity & Module & AClass & KSKernel & LClass & Run Time \\
\hline\hline
Threefold scroll & $S_{-2}$ & 1 & 0 & 1 & 0.2s \\
of type $(1,2)$ & $S_{-1}$ & 1 & 0 & 1 & 0.6s \\
& $S_1$ & 1 & 0 & 1 & 0.6s \\
& $M$ & 1 & 0 & 1 & 10.6s \\
\hline
Threefold Q.S. of & $S_{-1}$ & 1 & 0 & 0 & 8.3s \\
finite CM repr. type & $M$ & ? & 0 & ? & $\infty$ \\
\hline
Cubic cone & M3 & 1 & 0 & 0 & 0.9s \\
$x^3 + y^3 + z^3 = 0$ & M4 & 1 & 0 & 0 & 3.0s \\
\hline
Monomial curve & $p$ & 1 & 0 & 1 & 0.1s \\
\hline Threefold $E_6$ & $m$ & 1 & 0 & 0 & 0.4s \\
\hline
\end{tabular}
\end{center}

\bibliographystyle{elsart-harv}
\bibliography{main}

\end{document}